\documentclass [a4 paper,12pt]{article}
\usepackage{mathrsfs}

\usepackage[usenames]{color}
\usepackage{bbm,amsmath,amssymb,amsfonts,amsthm,graphics,graphicx}

\newtheorem{lem}{Lemma}
\newtheorem{thm}[lem]{Theorem}
\newtheorem{cor}[lem]{Corollary}

\title{The asymptotic number of occurrences of a subtree\\
in trees with bounded maximum degree\\ and an application to the
Estrada index\footnote{Supported by NSFC No.10831001, PCSIRT and the
``973" program.}}
\author{\small Xueliang Li, Yiyang Li\\
\small Center for Combinatorics and  LPMC-TJKLC\\
\small Nankai University, Tianjin 300071, China}
\date{ }

\begin{document}

\maketitle

\begin{abstract}
Let $\mathcal {T}^{\Delta}_n$ denote the set of trees of order $n$,
in which the degree of each vertex is bounded by some integer
$\Delta$. Suppose that every tree in $\mathcal {T}^{\Delta}_n$ is
equally likely. For any given subtree $H$, we show that the number
of occurrences of $H$ in trees of $\mathcal {T}^{\Delta}_n$ is with
mean $(\mu_H+o(1))n$ and variance $(\sigma_H+o(1))n$, where $\mu_H$,
$\sigma_H$ are some constants. As an application, we estimate the
value of the Estrada index $EE$ for almost all trees in $\mathcal
{T}^{\Delta}_n$, and give an explanation in theory to the
approximate linear correlation between $EE$ and the first Zagreb
index obtained by quantitative analysis.\\[3mm]
\\[3mm]
\end{abstract}

\section{Introduction}

We denote the set of trees with bounded maximum degree $\Delta$ by
$\mathcal {T}_n^\Delta$. Setting $t_n=|\mathcal{T}_n^\Delta|$, we
introduce a generating function for these trees
$$t(x)=\sum_{n\geq 1}t_n x^n.$$

Let $H$ be a given subtree. For a tree $T_n^\Delta$, we say that $H$
{\it occurs} in $T_n^\Delta$ if there is a subtree in $T_n^\Delta$
isomorphic to $H$. Denote the number of occurrence of $H$ in a tree
$T_n^\Delta$ by $t_n^\Delta$. To count the occurrences, we introduce
a generating function in two variables as follows
$$t(x,u)=\sum_{n\geq 1,T_n^\Delta\in \mathcal
{T}_n^\Delta}x^nu^{t_n^\Delta}.
$$
It can be simplified into
$$t(x,u)=\sum_{n\geq1, k\geq0}t_{n,k}x^nu^k,$$
where $t_{n,k}$ denotes the number of trees in
$\mathcal{T}_n^\Delta$ such that the number of occurrences of $H$ is
$k$. Note that $t(x,1)=t(x)$, i.e., $t_n=\sum_{k\geq0} t_{n,k}$.

Furthermore, suppose that every tree in $\mathcal{T}_n^\Delta$ is
equally likely. We use $X_n(T_n^\Delta)$, or simply $X_n$, to
denote the number of occurrences of $H$ in $T_n^\Delta$. Clearly,
$X_n$ is a random variable in $\mathcal{T}_n^\Delta$. Then, the
probability of $X_n$ can be defined as
$$\mbox{Pr}[X_n=k]=\frac{t_{n,k}}{t_n}.$$

Without the maximum degree restriction, many results have been
established for some special substructures. We further assume
that, when $H$ occurs in $T_n^\Delta$, the degrees of the internal
vertices (vertices of degrees not equal to 1) should coincide with
those of the corresponding vertices in the tree. In this case, $H$
is called a {\it pattern} of the tree. Kok \cite{gk} showed that
the number of occurrences of a pattern in the trees is with mean
$E (X_n)=(u+o(1))n$ and variance $Var (X_n)=(\sigma+o(1))n$, and
$\frac{X_n-E (X_n)}{\sqrt{Var (X_n)}}$ is asymptotic to a
distribution with density $(A+Bx^2)e^{-Cx^2}$ for some constants
$A, B, C\geq 0$. Moreover, if the pattern is a star, then the
number of occurrences of the pattern in a tree is exactly the
number of vertices with a given degree with respect to the
internal vertex of the star. It has been showed that for the
number of vertices of a given degree $X_n$, $\frac{X_n-E
(X_n)}{\sqrt{Var (X_n)}}$ is asymptotically normally distributed.
We refer the readers to \cite{dg, rs} for more details. And,
analogous results have been obtained for other classes of trees,
such as simply generated trees, rooted trees, {\it et al.} (see
\cite{cdkk}, \cite{dg} and \cite{gk}). However, for the number of
occurrences of a subtree $H$ in trees without the maximum degree
restriction, there is no such a similar result obtained. And, it
seems to be much difficult.

In this paper, we get that the number of occurrences of a subtree
$H$ in the planted trees and rooted trees with bounded maximum
degree is also asymptotically normally distributed with mean and
variance in $\Theta(n)$. And for $\mathcal {T}_n^\Delta$, we get a
weak result which does not show that the distribution is also
asymptotically normal. And, as an application, we use this result to
estimate the Estrada index $EE$ for these trees in $\mathcal
{T}_n^\Delta$, and give an explanation in theory to the approximate
linear correlation between $EE$ and the first Zagreb index
\cite{gfgmv} obtained by quantitative analysis.

In this paper, Section $2$ is devoted to a systematic treatment of
the number of occurrences of  a given subtree $H$. And in Section
$3$, we investigate the Estrada index for the trees in $\mathcal
{T}_n^\Delta$.

\section{The number of occurrences of a given subtree}

In this section, we proceed to show that the number of occurrences
of a subtree $H$ in $\mathcal {T}_n^\Delta$ is with mean
$(\mu_{H}+o(1))n$ and variance $(\sigma_H+o(1))n$ for some constants
$\mu_H$ and $\sigma_H$. And in the procedure of proof, we get the
related results for planted trees and rooted trees.

In what follows, we introduce some terminology and notation which
will be used in the sequel. For the others not defined here, we
refer to book \cite{hp}.

Analogous to trees, we introduce generating functions for rooted
trees and planted trees. Let $\mathcal{R}^{\Delta}_n$ denote the set
of rooted trees of order $n$ with degrees bounded by an integer
$\Delta$. Setting $r_n=|\mathcal {R}^{\Delta}_n|$, we have
$$r(x)=\sum_{n\geq1}r_n x^n$$ and $$r(x,u)=\sum_{n\geq1, k\geq0}r_{n,k}x^n
u^k,$$ where $r_{n,k}$ denotes the number of trees in
$\mathcal{R}^{\Delta}_n$ such that $H$ occurs $k$ times. A {\it
planted tree} is formed by adding a vertex to the root of a rooted
tree. The new vertex is called the {\it plant}, and we never count
it in the sequel. Analogously, let $\mathcal{P}^{\Delta}_n$ denote
the set of planted trees of order $n$ and bounded maximum degree
$\Delta$. Setting $p_n=|\mathcal{P}^{\Delta}_n|$, we have
$$p(x)=\sum_{n\geq1}p_n x^n$$ and $$p(x,u)=\sum_{n\geq 1,
k\geq0}p_{n,k}x^n u^k,$$ where $p_{n,k}$ denotes the number of trees
in $\mathcal{P}^{\Delta}_n$ such that $H$ occurs $k$ times. By the
definition of planted trees, one can readily see that
$p(x,1)=p(x)=r(x,1)=r(x)$.

Moreover, in \cite{ot}, it has been showed that there exists a
number $x_0$ such that
\begin{equation}\label{expand}p(x)=b_1+b_2\sqrt{x_0-x}+b_3
(x_0-x)+\cdots,\end{equation} where $b_1,b_2,b_3$ are some constants
not equal to zero. Evidently, $p(x_0)=b_1$ and for any $|x|\leq
x_0$, $p(x)$ is convergent. For any $\Delta$, $x_0\leq1/2$,
particularly, if $\Delta=4$, $x_0\approx 0.3551817$ and
$p(x_0)\approx1.117421$.

Let $p^{(\Delta-1)}(x)$ be the generating function of trees such
that the degrees of the roots are not more than $\Delta-1$, while
the degrees of the other vertices are still bounded by $\Delta$.
Then, we have (see \cite{ot})
\begin{equation}\label{De-1}
p^{(\Delta-1)}(x_0)=1.
\end{equation} And, this fact will play an important role in the
following proof.

Let $\textbf{y}(x,u)=(y_1(x,u),\ldots,y_N(x,u))^T$ be a column
vector. We suppose that $G(x,\textbf{y},u)$ is an analytic function
with non-negative  Taylor coefficients. $G(x,\textbf{y},u)$ can be
expanded as
$$ G(x,\textbf{y},u)=\sum_{n\geq1,k\geq0}g_{n,k}x^nu^k.$$
Let $X_n$ denote a random variable with probability
\begin{equation}\label{rv}
\mbox{Pr}[X_n=k]=\frac{g_{n,k}}{g_n},
\end{equation}
where $g_n=\sum_kg_{n,k}$. First, we introduce a useful lemma
\cite{cdkk, d}.

\begin{lem}\label{mainlem}
Let
$\textbf{F}(x,\textbf{y},u)=(F^1(x,\textbf{y},u),\ldots,F^N(x,\textbf{y},u))^T$
be functions analytic around $x=0$,
$\textbf{y}=(y_1,\ldots,y_N)^T=\textbf{0}$, $u=0$, with Taylor
coefficients all are non-negative. Suppose
$\textbf{F}(0,\textbf{y},u)=\textbf{0}$,
$\textbf{F}(x,\textbf{0},u)\neq\textbf{0}$,
$\textbf{F}_x(x,\textbf{y},u)\neq \textbf{0}$, and for some $j$,
$\textbf{F}_{y_jy_j}(x,\textbf{y},u)\neq\textbf{0}$. Furthermore,
assume that $x=x_0$ together with $\textbf{y}=\textbf{y}_0$ is a
non-negative solution of the system of equations
\begin{align}
&\textbf{y}=\textbf{F}(x,\textbf{y},1)\\
&0=\mbox{det}(\textbf{I}-\textbf{F}_\textbf{y}(x,\textbf{y},1))\label{det}
\end{align}
inside the region of convergence of $\textbf{F}$, $\textbf{I}$ is
the unit matrix. Let $\textbf{y}=(y_1(x,u),\ldots,y_N(x,u))^T$
denote the analytic solution of the system
\begin{equation}\label{fs}
\textbf{y}=\textbf{F}(x,\textbf{y},u)
\end{equation}
with $\textbf{y}(0,u)=\textbf{0}$.

If the dependency graph $G_{\textbf{F}}$ of the function system
Equ.(\ref{fs}) is strongly connected, then there exist functions
$f(u)$ and $g_i(x,u)$, $h_i(x,u)$ ($1\leq i\leq N$) which are
analytic around $x=x_0$, $u=1$, such that
\begin{equation}\label{sqrt}
y_i(x,u)=g_i(x,u)-h_i(x,u)\sqrt{1-\frac{x}{f(u)}}
\end{equation}
is analytically continued around $u=1$, $x=f(u)$ with
$\mbox{arg}(x-f(u))\neq 0$, where $x=f(u)$ together with
$y=y(f(u),u)$ is the solution of the extended system
\begin{align}
&\textbf{y}=\textbf{F}(x,\textbf{y},u)\\
&0=\mbox{det}(\textbf{I}-\textbf{F}_\textbf{y}(x,\textbf{y},u))\label{det2}.
\end{align}

Moreover, let $G(x,\textbf{y},u)$ be an analytic function with
non-negative Taylor coefficients such that the point
$(x_0,\textbf{y}(x_0,1),1)$ is contained in the region of
convergence. Finally, let $X_n$ be the random variable defined in
Equ.(\ref{rv}). Then the random variable $X_n$ is asymptotically
normal with mean
$$E(X_n)=\mu n+O(1)\mbox{  } (n\rightarrow \infty),$$ and variance
$$Var(X_n)=\sigma n+O(1)\mbox{ } (n\rightarrow \infty)$$
with $\mu=\frac{-f'(1)}{f(1)}$.
\end{lem}

\noindent {\bf Remark 1:} We say that the {\it dependency graph}
$G_{\textbf{F}}$ of $\textbf{y}=\textbf{F}(x,\textbf{y},u)$ is
strongly connected if there is no subsystem of equations that can be
solved independently from others. If $G_\textbf{F}$ is strongly
connected, then
$\textbf{I}-\textbf{F}_\textbf{y}(x_0,\textbf{y}_0,1)$ has rank
$N-1$. Suppose that $\mathbf{v}^T$ is a vector with
$\mathbf{v}^T(\mathbf{I}-\mathbf{F}_\mathbf{y}(x_0,\mathbf{y_0},1))=0$.
Then, $\mu=\frac{\mathbf{v}^T(\mathbf{F}_u(x_0,\mathbf{y_0},1))}
{x_0\mathbf{v}^T(\mathbf{F}_x(x_0,\mathbf{y_0},1))}$. We refer the
readers to \cite{cdkk, d} for more details.

Now, we concentrate on considering the generating function $p(x,u)$.

For the subtree $H$, we suppose the diameter of $H$ is $h$. The {\it
depth} of a vertex in a planted tree is the distance from the vertex
to the root. The {\it depth of a planted tree} is the largest
distance from the vertices to the root. We split up
$\mathcal{P}_n^{\Delta}$ into two sets $\mathcal{W}_0$ and
$\mathcal{W}$, which denote the trees with depth not more than $h-1$
and the trees with depth greater than $h-1$, respectively. We can
see that if $H$ occurs in the planted tree and the corresponding
subtree takes the root, then the depth of the subtree is not more
than $h$. Moreover, since we mainly consider the asymptotic number
of subtree, the trees in $\mathcal{W}_0$ will contribute nothing to
the coefficient of $x^nu^k$ when $n$ is large enough. Therefore, in
this paper, we do not need to know the exactly expression of the
generating function for the trees in $\mathcal{W}_0$, and denote it
by $\phi(x,u)$. Moreover, the counting function of a subset in
$\mathcal{W}_0$ is also denoted by $\phi(x,u)$. In what follows, we
shall see that this assumption is reasonable. Then, we focus on the
trees in $\mathcal{W}$.

First, we introduce some conceptions. For a planted tree in
$\mathcal{W}$, the planted subtree formed by the vertices with depth
not more than $\ell$ is called {\it $\ell$-depth subtree} of this
tree. Now, we split up $\mathcal {W}$ according to the $h$-depth
subtree. That is, the trees in $\mathcal {W}$ having the same
$h$-depth subtree $w_i$ form a subset $\mathcal {H}_i$ of $\mathcal
{W}$. Since the degrees of the vertices in $\mathcal {W}$ are
bounded by $\Delta$, there are finite number of different $h$-depth
subtrees, i.e., for some $N_\Delta$, $1\leq i\leq N_\Delta$.
Therefore, we obtain that
\begin{equation}\label{p(x,u)}p(x,u)=\phi(x,u)+\sum_{i=1}^{N_\Delta}a_{w_i,h}(x,u),\end{equation}
where $a_{w_i,h}(x,u)$ denotes the generating function of $\mathcal
{H}_i$.

To establish the functions system of $a_{w_i,h}(x,u)$, we need other
functions $a_{w'_i,h-1}(x,u)$ as follows. For some $(h-1)$-depth
subtree $w'_i$, we denote the subset of the planted trees in
$\mathcal{W}$ having $w'_i$ by $\mathcal{H}_i'$. Note that
$w'_i\notin \mathcal{H}_i'$. Then, the set of the planted trees
having $w'_i$ consists of $\mathcal {H}'_i$ and $w'_i$. If we use
$a_{w'_i,h-1}(x,u)$ to denote the generating function with respect
to $\mathcal{H}_i'\cup {w'_i}$, it follows that
\begin{equation}\label{prefunction}a_{w'_i,h-1}(x,u)=\sum_{w_i\in \mathcal{H}_i'}
a_{w_i,h}(x,u)+w'_i(x,u).\end{equation}

There will appear an expression of the form $Z(S_n, f(x,u))$ (or
$f(x)$), which is the substitution of the counting series $f(x,u)$
(or $f(x)$) into the cycle index $Z(S_n)$ of the symmetric group
$S_n$. This involves replacing each variable $s_i$ in $Z(S_n)$ by
$f(x^i,u^i)$ (or $f(x^i)$). For instance, if $n=3$, then
$Z(S_3)=(1/3!)(s_1^3+3s_1s_2+2s_3)$, and
$Z(S_3,f(x,u))=(1/3!)(f(x,u)^3+3f(x,u)f(x^2,u^2)+2f(x^3,u^3))$. We
refer the readers to \cite{hp} for details.

Note that a planted tree can be seen as a root attached by some
planted subtrees. Employing the classic P\'{o}lya enumeration
theorem, we have $Z(S_{j-1};p(x))$ as the counting series of the
planted trees whose roots have degree $j$, and the coefficient of
$x^p$ in $x\cdot Z(S_{j-1};p(x))$ is the number of planted trees
with $p$ vertices (see \cite{hp} p.51--54). Therefore,
$$p(x)=x\cdot\sum_{j=0}^{\Delta-1}Z(S_j;p(x)),$$ and
$$p^{(\Delta-1)}(x)=x\cdot\sum_{j=0}^{\Delta-2}Z(S_j;p(x)).$$

By means of the same method, $a_{w_i,h}(x,u)$ can be expressed in
$a_{w'_i,h-1}(x,u)$ and $\phi(x,u)$. Suppose that the roots of the
trees in $\mathcal{H}_i$ have degree $j$, and each has $j'$ planted
subtrees with depth at least $h-1$ attached to it. Clearly, $j'$
belongs to $\{1,\ldots,j-1\}$, and some of these subtrees may have
the same $w'_i$. Denote these different $(h-1)$-depth subtrees by
$\{w'_{s}\}$ and suppose $w'_{s}$ happens $\ell_s$ times. Evidently,
$\sum \ell_s=j'$. It follows that
\begin{equation}\label{functions}a_{w_i,h}(x,u)=x\cdot
\prod_{} Z(\ell_s; a_{w'_s,h-1})\cdot \phi(x,u)\cdot u^{k(\ell_s,
\phi)} \mbox{, }(1\leq i\leq N_\Delta).\end{equation} Here,
$\phi(x,u)$ denotes the counting function of the other $j-1-j'$
planted subtrees, since these subtrees belong to $\mathcal{W}_0$.
And the factor $u^{k(\ell_s, \phi)}$ serves to count the number of
occurrences of $H$ using the root of the new tree. In this case, all
these vertices of the new tree corresponding the vertices of $H$
have depth not more than $h$. And, since we know that the $h$-depth
subtree of the new tree is $w_i$, the number of occurrences taking
the root can be calculated, that is, the upper index $k(l_s,\phi)$
can be calculated. Therefore, combining with
Equ.(\ref{prefunction}), the functions system of $a_{w_i,h}(x,u)$
has been established.

Now, we start to show that all the conditions of Lemma \ref{mainlem}
hold for $a_{w_i,h}(x,u)$. For convenience, we still use
$\mathbf{F}$ to denote the functions system. Set vector
$\mathbf{a}(x,u)=(a_{w_1,h},\ldots,a_{w_{N_\Delta},h})$. We suppose
$a_{w_i,h}(x,u)=F^i(x,\mathbf{a},u)$. Since $p(x,1)=p(x)$ and
$p(x_0)=b_1$, one can see that $a_{w_i,h}(x_0,1)$ is convergent. So,
$x_0$ and $\mathbf{a}(x_0,1)$ are inside the region of convergence
of $\mathbf{F}$. Apparently, the other conditions are easy to verify
except for Equ.(\ref{det}). In what follows, we shall show that the
sum $S_{a_{w_i,h}}$ of every column of $\mathbf{F}_\mathbf{a}(x_0,
\mathbf{a}(x_0,1),1)$ equals $1$. Consequently, the equation
$\mbox{det}(\mathbf{I}-\mathbf{F}_\mathbf{a}(x_0,
\mathbf{a}(x_0,1),1))=0$ holds.

We consider the derivative on $a_{w_{i_0},h}$. If
$F^i(x,\mathbf{a},u)$ is not the function of $a_{w_{i_0},h}$, then
$F^i_{a_{w_{i_0},h}}(x,\mathbf{a},u)$ will contribute nothing to the
sum $S_{a_{w_i,h}}$. Thus, we just need to consider the functions
$F^i(x,\mathbf{a},u)$ with some $a_{w'_s,h-1}$ having the term
$a_{w_{i_0},h}$. In Equ.({\ref{functions}), if both
$a_{w'_{s_1},h-1}$ and $a_{w'_{s_2},h-1}$ have the term
$a_{w_{i_0},h}$, which implies that the trees corresponding to
$a_{w'_{s_1},h-1}$, $a_{w'_{s_2},h-1}$ have the same $(h-1)$-depth
subtree, then by the definition of $a_{w'_{s},h-1}$, we get that
$a_{w'_{s_1},h-1}=a_{w'_{s_2},h-1}$. Therefore, there exists exactly
one product factor, say $Z(l_{s_0}; a_{w'_{s_0},h-1})$, in the
expression of $a_{w_{i_0},h}$.

Moreover, it is well-known that the partial derivative of
$Z(S_n;\cdot)$ enjoys (see \cite{dg})
\begin{equation}\label{pdz}
\frac{\partial}{\partial s_1}Z(S_n;
s_1,\ldots,s_n)=Z(S_{n-1};s_1,\ldots,s_{n-1}).
\end{equation}
For the planted tree, we have $\frac{\partial Z(S_n;
p(x,1))}{\partial p(x,1)}=Z(S_{n-1}; p(x,1))$, which equals the
generating function obtained by deleting one subtree from the root.
Analogously, we have
$$F^i_{a_{w_{i_0},h}}=x\cdot \prod_{s\neq s_0} Z(\ell_s;
a_{w'_s,h-1})\cdot Z(\ell_{s_0}-1;
a_{w'_{s_0},h-1})\cdot\phi(x,u)\cdot u^{k(\ell_s, \phi)}, $$ and it
is exactly the new generating function produced by deleting one
planted subtree of $\mathcal{H}'_{s_0}\cup w'_{s_0}$. Clearly, the
root of the new planted tree is of degree $j-1$. Particularly, if
$\ell_{s_0}=1$, after taking the derivative, the yielded function
corresponds to the trees with roots of degree $j-1$ such that every
planted subtree does not belong to $\mathcal {H}'_{s_0}\cup
w'_{s_0}$. Set $u=1$. It follows that $S_{a_{w_i,h}}(x,
\mathbf{a}(x,1),1)$ is the generating function of the planted trees
with roots of degree not more than $\Delta-1$, i.e.,
$p^{(\Delta-1)}(x,1)$. Combining with the fact
$p^{(\Delta-1)}(x_0,1)=1$, we obtain $S_{a_{w_i,h}}(x_0,
\mathbf{a}(x_0,1),1)=1$. Immediately, the Equ.(\ref{det})
$$\mbox{det}(\mathbf{I}-\mathbf{F}_\mathbf{a}(x_0,
\mathbf{a}(x_0,1),1))=0$$ follows.

Employing Lemma \ref{mainlem}, we have that $a_{w_i,h}(x,u)$ is in
the form of Equ.(\ref{sqrt}), namely, for some $f(u)$ and
$g_{w_i,h}(x,u)$, $h_{w_i,h}(x,u)$ which are analytic around
$x=x_0$, $u=1$, it follows that
\begin{equation*}
a_{w_i,h}(x,u)=g_{w_i,h}(x,u)-h_{w_i,h}(x,u)\sqrt{1-\frac{x}{f(u)}}
\end{equation*} is analytically continued around $u=1$, $x=f(u)$ with
$\mbox{arg}(x-f(u))\neq 0$. From Equ.(\ref{p(x,u)}), we can see that
$p(x,u)$ can be written into a function of $\mathbf{a}(x,u)$, and
denote it by $P(x,\mathbf{a}(x,u),u)$. Clearly, all the coefficients
of $P(x,\mathbf{a}(x,u),u)$ are non-negative. Therefore, $p(x,u)$ is
also in the form of Equ.(\ref{sqrt}). Moreover, recalling
Equ.(\ref{expand}), we can see that $f(1)=x_0$. Apply Lemma
\ref{mainlem} to $P(x,\mathbf{a}(x,u),u)$, the following result is
obtained.

\begin{thm}
For any given subtree $H$, the number $X_n$ of occurrences of $H$ in
$\mathcal{P}_n^\Delta$ is asymptotical to be normal with mean
$E(X_n)=\mu_H n+O(1)$ and variance $Var (X_n)=\sigma_H^p n+O(1)$ for
some constants $\mu_H$ and $\sigma_H^p$.
\end{thm}

A rooted tree in $\mathcal {R}_n^{\Delta}$ can also be seen as a
root attached by some planted trees. That is, by the classic
P\'{o}lya enumeration theorem, analogous to Equ.(\ref{functions}),
the generating function of $\mathcal {R}_n^{\Delta}$ is also a
function in $\mathbf{a}(x,u)$. We denote the function by
$R(x,\mathbf{a}(x,u),u)$, and $r(x,u)=R(x,\mathbf{a}(x,u),u)$. By
means of the above analysis, it is not difficult to see that the
Taylor coefficients of $R(x,\mathbf{a}(x,u),u)$ are non-negative.
Thus, $r(x,u)$ also has the form of Equ.(\ref{sqrt}). And, apply
Lemma \ref{mainlem} to $R(x,\mathbf{a}(x,u),u)$, the following
result is obtained.
\begin{thm}
For any given subtree $H$, the number $X_n$ of occurrences of $H$ in
$\mathcal{R}_n^\Delta$ is asymptotically normally distributed with
mean $E(X_n)=\mu_H n+O(1)$ and variance $Var(X_n)=\sigma_H^r n+O(1)$
for some constants $\mu_H$ and $\sigma_H^r$.
\end{thm}

\noindent{\bf Remark 2:}  Since $r(x,u)$ and $p(x,u)$ correspond
to the same function $f(u)$, by Lemma \ref{mainlem} we can see
that the means of $X_n$ with respect to $\mathcal{R}_n^\Delta$ and
$\mathcal{P}_n^\Delta$ are with the same constant $\mu_H$.
Moreover, it has been showed that the sum of each column of
$\mathbf{F}_{\mathbf{a}}(x_0,\mathbf{a}(x_0,1),1)$ equals $1$,
then we have $\mathbf{v}^T=(1,\ldots,1)$ such that
$\mathbf{v}^T(\mathbf{I}-\mathbf{F}_\mathbf{y}(x_0,\mathbf{y_0},1))=0$.
Therefore, it is easy to see that $\mu_H$ is positive by Remark $1$. \\

In what follows, we investigate the generating function of trees.
Two edges in a tree are {\it similar}, if they are the same under
some automorphism of the tree. To {\it join} two planted trees is to
connect the two roots with a new edge and get rid of the two plants.
If the two panted trees are the same, we say that the new edge is
{\it symmetric}. Then, we have the following lemma due to \cite{ot}.
\begin{lem}\label{otter}
For any tree, the number of rooted trees corresponding to this tree
minus the number of nonsimilar edges (except for the symmetric edge)
is the number $1$.
\end{lem}

Note that, if we delete any one edge from a similar set in a tree,
the yielded trees are the same two trees. Hence, different pairs of
planted trees correspond to nonsimilar edges. Now, we have
\begin{align}\label{t(x,u)}
t(x,u)=&r(x,u)-\frac{1}{2}\Big(\sum_{1\leq i_1,i_2\leq N_{\Delta}}
a_{w_{i_1},h}(x,u)a_{w_{i_2},h}(x,u)\cdot
u^{k(w_{i_1},w_{i_2})}\Big)\nonumber\\
&+\frac{1}{2}\sum_{1\leq i\leq N_{\Delta}}a_{w_{i},h}(x^2,u^2)\cdot
u^{k(w_i,w_i)},
\end{align}
where $k(w_{i_1},w_{i_2})$ serves to count the subtrees taking
vertices both in $w_{i_1}$ and $w_{i_2}$. Consequently, we obtain
that $t(x,u)$ is also in the form of Equ.(\ref{sqrt}), i.e., there
exist some functions $\overline{g}(x,u)$, $\overline{h}(x,u)$ which
are analytic around $x=x_0$, $u=1$, such that
\begin{align*}
t(x,u)&=\overline{g}(x,u)-\overline{h}(x,u)\sqrt{1-\frac{x}{f(u)}}.
\end{align*}
is analytically continued around $u=1$, $x=f(u)$ with
$\mbox{arg}(x-f(u))\neq 0$. Here, we could not show that $t(x,u)$
likes $P(x,a(x,u),u)$ and $R(x,a(x,u),u)$ that have non-negative
Taylor coefficients, so Lemma \ref{mainlem} fails in this case.
However, we can use the following result due to \cite{gk} to get a
weak result for $t(x,u)$.
\begin{lem}\label{lem3}
Suppose that $t(x,u)$ has the form
$$t(x,u)=\overline{g}(x,u)-\overline{h}(x,u)\sqrt{1-\frac{x}{f(u)}},$$
where $\overline{g}(x,u)$, $\overline{h}(x,u)$ and $f(u)$ are
analytic functions around $x=f(1)$ and $u=1$ that satisfy
$\overline{h}(f(1),1)=0$, $\overline{h}_x(f(1),1)\neq 0$, $f(1)>0$
and $f'(1)<0$. Furthermore, $x=f(u)$ is the only singularity on the
cycle $|x|=|f(u)|$ for $u$ is close to $1$. Suppose that $X_n$ is
defined as Equ.(\ref{rv}) to $y(x,u)$. Then, $E(X_n)=(\mu+o(1))n$
and $Var(X_n)=(\sigma+o(1))n$, where $\mu=-f'(1)/f(1)$ and
$\sigma=\mu^2+\mu-f''(1)/f(1)$.
\end{lem}

\noindent {\bf Remark 3:} This result does not tell us that the
limiting distribution is asymptotically normal. If $h(f(1),1)\neq
0$, this lemma is trivial by Lemma \ref{mainlem}, and if
$h(f(u),u)=0$, we can still get that the limiting distribution is
normal by further analysis (see \cite{dg}).\\

For $t(x)$, it has been obtained that \cite{ot}
$$t(x)=c_0+c_1(x_0-x)+c_2(x_0-x)^{3/2}+\cdots, $$
where $c_0$, $c_1$, $c_2$  are some constants not equal to $0$.
Combining with the fact $t(x,1)=t(x)$, we can see that
$\overline{h}(f(1),1)=0$ and $\overline{h}_x(f(1),1)\neq0$.
Moreover, the other conditions in Lemma \ref{lem3} are easy to
verify. Then, we formulate the following theorem.

\begin{thm}
Let $X_n$ be the number of occurrences of a given subtree $H$ in
the trees of $\mathcal{T}_n^{\Delta}$. Then it follows that
$$E(X_n)=(\mu_H+o(1)) n$$ and
$$Var(X_n)=(\sigma_H^t+o(1))n,$$ where $\mu_H$ and
$\sigma_H^t$ are some constants with respect to the subtree $H$.
\end{thm}

Following book \cite{BB}, we will say that {\it almost every} (a.e.)
graph in a random graph space $\mathcal{G}_n$ has a certain property
$Q$ if the probability $\mbox{Pr}(Q)$ in $\mathcal{G}_n$ converges
to 1 as $n$ tends to infinity. Occasionally, we shall write {\it
almost all} instead of almost every.

By Chebyshev inequality one can get that
$$ \mbox{Pr}\big[\big|X_n-E(X_n)\big|>n^{3/4}\big]\leq \frac{Var
(X_n)}{n^{3/2}}\rightarrow 0 \mbox{ as } n\rightarrow \infty.$$
Therefore, for any subtree $H$, $X_n=(\mu_H+o(1)) n \mbox{  a.e.} $
in $\mathcal{T}_n^\Delta$. Then, an immediate consequence is the
following.

\begin{cor}
For almost all trees in $\mathcal{T}_n^\Delta$, the number of
occurrences of $H$ equals $(\mu_H+o(1)) n$.
\end{cor}

\section{The Estrada index}

In this section, we explore the Estrada index for trees in
$\mathcal{T}_n^{\Delta}$. Let $G$ be a simple graph with $n$
vertices. The eigenvalues of the adjacency matrix of $G$ are said to
be the eigenvalues of $G$ and to form the spectrum. Suppose that the
eigenvalues of $G$ are $\lambda_i$, $1\leq i\leq n$. The Estrada
index is defined as
$$EE(G)=\sum_{i=1}^n e^{\lambda_i}.$$
This index is invented in year $2000$, and nowadays widely accepted
and used in the information-theoretical and network-theoretical
applications. And for this graph invariant, many results have been
established. We refer the readers to a survey \cite{drg} for more
details.

Furthermore, for trees with $n$ vertices, it has been showed that
the path has the minimum Estrada index and the star has the maximum.
And by quantitative analysis, there is an approximate linear
correlation between $EE$ and the first Zagreb index, i.e., $\sum
d_i^2$ for trees. Denote $\sum d_i^2$ by $D$. That is,
\begin{equation}\label{ee}EE\approx a D+b,\end{equation} where $a$
and $b$ are some constants. We refer the readers to \cite{drg} and
\cite{gfgmv}.

In what follows, we shall get the estimate of $EE$ for almost all
trees in $\mathcal{T}_n^\Delta$ and give an explanation to the
correlation (\ref{ee}) in theory.

Denoting by $M_k=M_k(G)=\sum_{i=1}^{n}\lambda_i^k$ the $k$-th
spectral moment of $G$, and bearing in mind the power-series
expansion of $e^x$, we have
$$EE(G)=\sum_{k=0}^{\infty}\frac{M_k(G)}{k!}.$$
Note that $M_k(G)$
is equal to the number of closed walks of length $k$. For trees, one
can readily see that
\begin{equation}\label{eet}EE(T)=\sum_{k=0}^{\infty}\frac{M_{2k}}{(2k)!}.\end{equation}

Then, in a tree, the closed walk of length $2k$ forms a subtree with
at most $k+1$ vertices. We have got that, for any given subtree, the
number of occurrences of the subtree in $\mathcal{T}_n^{\Delta}$
equals $(\mu_{H}+o(1))n \mbox{ a.e.}$ Since there are finite
different subtrees with at most $k+1$ vertices, and each subtree
corresponds to finite numbers of $2k$ closed walks, we can obtain
that there exists a constant $\mu_{2k}$ such that the number of $2k$
closed walk is $(\mu_{2k}+o(1))n \mbox{ a.e.}$, namely,
$$M_{2k}=(\mu_{2k}+o(1))n \mbox{ a.e.}$$ in
$\mathcal{T}_n^{\Delta}$. Moreover, we introduce a lemma due to Fiol
and Garriga \cite{fg}.
\begin{lem}\label{walkdegree}
For any graph $G$,  $M_{2k}\leq \sum_{i=1}^n d_i^{2k}$
\end{lem}

Recall that the degrees of a tree in $\mathcal{T}_n^{\Delta}$ are
bounded by $\Delta$. So, $\sum_{i=1}^n d_i^{2k}\leq \Delta^{2k}n$
and thus $EE(T_n^\Delta)\leq e^\Delta n$. Moreover, since
$\sum_{k=0}\frac{\Delta ^{2k}}{(2k)!}$ is convergent, for any
positive number $\varepsilon$, there exists an integer $j_0$ such
that for any $j>j_0$, $\sum_{k=j+1}\frac{M_{2k}}{(2k)!}<\varepsilon
n$. Evidently, it is uniform for all the trees in
$\mathcal{T}_n^{\Delta}$. Therefore, we have
$$\sum_{k=0}^{j}\frac{M_{2k}}{(2k)!}\leq EE(T_n^{\Delta})\leq
\sum_{k=0}^{j}\frac{M_{2k}}{(2k)!}+\varepsilon n.$$ Hence, we just
contribute to consider the closed walks of length at most $j_0$.

For any integer $j$, we have
$\sum_{k=0}^{j}\frac{\mu_{2k}}{(2k)!}\leq e^{\Delta}$. Therefore,
$\sum_{k=0}\frac{\mu_{2k}}{(2k)!}$ is convergent, and denote the
limit by $\mu_\Delta$. It follows that
$$(\mu_{\Delta}-\varepsilon)n<\sum_{k=0}^{j}\frac{M_{2k}}{(2k)!}
=\sum_{k=0}^{j}\frac{(\mu_{2k}+o(1))n}{(2k)!}\leq (\mu_\Delta+o(1))
n \mbox{ a.e.}$$ Then, we have that
$(\mu_\Delta-\varepsilon)n<EE(T_n^{\Delta})<(\mu_\Delta+\varepsilon)n
\mbox{ a.e.}$ Now, we can formulate the following theorem.

\begin{thm}\label{correlation}
For any $\varepsilon>0$, the Estrada index of a tree in $\mathcal
{T}_n^{\Delta}$ enjoys
$$(\mu_\Delta-\varepsilon)n<EE(T_n^{\Delta})<(\mu_\Delta+\varepsilon)n \mbox{ a.e.},$$
where $\mu_{\Delta}$ is some constant.
\end{thm}

If we suppose that the given subtree $H$ is a path $L$ of length
$2$, then there exists some constant $u_L$ such that in
$\mathcal{T}_n^{\Delta}$, the number of occurrences $X_n$ of $L$ is
$(u_L+o(1))n\mbox { a.e.}$ In this case, it is easy to see that for
each tree $T_n^{\Delta}$, $X_n(T_n^{\Delta})=\sum_i {d_i\choose
2}=\frac{1}{2} D(T_n^{\Delta})-n+1$. Therefore, the value of $D$
also enjoys $(u_D+o(1))n\mbox { a.e.}$ for some constant $u_D$.
Then, combining with Theorem \ref{correlation}, we can see that, for
trees in $\mathcal{T}_n^\Delta$, the correlation between $EE$ and
$D$ is approximate to be linear.

\end{document}